

\input harvmac.tex
\input epsf.tex

\overfullrule=0mm
\hfuzz 15pt

\newcount\figno
\figno=0
\def\fig#1#2#3{
\par\begingroup\parindent=0pt\leftskip=1cm\rightskip=1cm\parindent=0pt
\baselineskip=11pt
\global\advance\figno by 1
\midinsert
\epsfxsize=#3
\centerline{\epsfbox{#2}}
\vskip 12pt
{\bf Fig. \the\figno:} #1\par
\endinsert\endgroup\par
}
\def\figlabel#1{\xdef#1{\the\figno}}
\def\encadremath#1{\vbox{\hrule\hbox{\vrule\kern8pt\vbox{\kern8pt
\hbox{$\displaystyle #1$}\kern8pt}
\kern8pt\vrule}\hrule}}


%

%
%
\def\frac#1#2{{\scriptstyle{#1 \over #2}}}

\def\ket#1{ | #1 \rangle}
\def\bra#1{ \langle #1 |}
%
%
\def\CA{{\cal A}}              \def\CC{{\cal C}}
              
\def\CG{{\cal G}}       \def\CH{{\cal H}}

\def\CP{{\cal P}}              
              \def\CU{{\cal U}}

\def\({ \left( }\def\[{ \left[ }
\def\){ \right) }\def\]{ \right] }
%


\def\IR{\relax{\rm I\kern-.18em R}}
\font\cmss=cmss10 \font\cmsss=cmss10 at 7pt
\def\IZ{\relax\ifmmode\mathchoice
{\hbox{\cmss Z\kern-.4em Z}}{\hbox{\cmss Z\kern-.4em Z}}
{\lower.9pt\hbox{\cmsss Z\kern-.4em Z}}
{\lower1.2pt\hbox{\cmsss Z\kern-.4em Z}}\else{\cmss Z\kern-.4em Z}\fi}
\def\inbar{\,\vrule height1.5ex width.4pt depth0pt}
\def\IB{\relax{\rm I\kern-.18em B}}
\def\ID{\relax{\rm I\kern-.18em D}}
\def\IE{\relax{\rm I\kern-.18em E}}
\def\IF{\relax{\rm I\kern-.18em F}}
\def\IG{\relax\hbox{$\inbar\kern-.3em{\rm G}$}}
\def\IH{\relax{\rm I\kern-.18em H}}
\def\II{\relax{\rm I\kern-.18em I}}
\def\IK{\relax{\rm I\kern-.18em K}}
\def\IL{\relax{\rm I\kern-.18em L}}
\def\IM{\relax{\rm I\kern-.18em M}}
\def\IN{\relax{\rm I\kern-.18em N}}
\def\IO{\relax\hbox{$\inbar\kern-.3em{\rm O}$}}
\def\IP{\relax{\rm I\kern-.18em P}}
\def\IQ{\relax\hbox{$\inbar\kern-.3em{\rm Q}$}}
\def\IGa{\relax\hbox{${\rm I}\kern-.18em\Gamma$}}
\def\IPi{\relax\hbox{${\rm I}\kern-.18em\Pi$}}
\def\ITh{\relax\hbox{$\inbar\kern-.3em\Theta$}}
\def\IOm{\relax\hbox{$\inbar\kern-3.00pt\Omega$}}


\def\oh{{1\over 2}}

\def\Ga{\alpha}\def\Gb{\beta}\def\Gc{\gamma}
\def\Gd{\delta}\def\GD{\Delta}\def\Ge{\epsilon}

\def\Gl{\lambda}\def\GL{\Lambda}
\def\Gm{\mu}\def\Gn{\nu}

\def\Gs{\sigma}


\def\mod{{\rm mod\,}}

\def\bra{\langle}\def\ket{\rangle}
\def\nind{\noindent}

\def\hepth#1{{\tt hep-th #1}}

\def\\#1 {{\tt\char'134#1} }

\catcode`\@=11
\def\Eqalign#1{\null\,\vcenter{\openup\jot\m@th\ialign{
\strut\hfil$\displaystyle{##}$&$\displaystyle{{}##}$\hfil
&&\qquad\strut\hfil$\displaystyle{##}$&$\displaystyle{{}##}$
\hfil\crcr#1\crcr}}\,}   \catcode`\@=12

\def\encadre#1{\vbox{\hrule\hbox{\vrule\kern8pt\vbox{\kern8pt#1\kern8pt}
\kern8pt\vrule}\hrule}}
\def\encadremath#1{\vbox{\hrule\hbox{\vrule\kern8pt\vbox{\kern8pt
\hbox{$\displaystyle #1$}\kern8pt}
\kern8pt\vrule}\hrule}}

\def\bun{{\bf 1}}


\def\IC{\relax\hbox{$\inbar\kern-.3em{\rm C}$}}

\def\msy{y }
\ifx\msan\msy
\input amssym.def
\input amssym.tex
\def\IZ{\Bbb Z}\def\IR{\Bbb R}\def\IC{\Bbb C}\def\IN{\Bbb N}

\else

\fi


\def\qm{q^{-1}}
\def\tr{{\rm tr\,}}
\def\Zt{\widetilde{Z}}
\def\l{{\ell}}

\newdimen\xraise\newcount\nraise
\def\xpoint{\hbox{\vrule height .45pt width .45pt}}
\def\udiag#1{\vcenter{\hbox{\hskip.05pt\nraise=0\xraise=0pt
\loop\ifnum\nraise<#1\hskip-.05pt\raise\xraise\xpoint
\advance\nraise by 1\advance\xraise by .4pt\repeat}}}
\def\ddiag#1{\vcenter{\hbox{\hskip.05pt\nraise=0\xraise=0pt
\loop\ifnum\nraise<#1\hskip-.05pt\raise\xraise\xpoint
\advance\nraise by 1\advance\xraise by -.4pt\repeat}}}
\def\bdiamond#1#2#3#4{\raise1pt\hbox{$\scriptstyle#2$}
\,\vcenter{\vbox{\baselineskip12pt
\lineskip1pt\lineskiplimit0pt\hbox{\hskip10pt$\scriptstyle#3$}
\hbox{$\udiag{30}\ddiag{30}$}\vskip-1pt\hbox{$\ddiag{30}\udiag{30}$}
\hbox{\hskip10pt$\scriptstyle#1$}}}\,\raise1pt\hbox{$\scriptstyle#4$}}



\lref\Res{N.Yu. Reshetikhin,
{\it Quantized Universal Enveloping Algebras, the Yang-Baxter 
Equation and Invariants of Links}, I and II, LOMI preprints 1988. }

\lref\BaRe{V.V. Bazhanov and N.Yu. Reshetikhin, 
{\it J. Phys.}{\bf A 23} (1990) 1477.}

\lref\VPun{V. Pasquier, {\it Nucl. Phys.} {\bf B285} [FS19] (1987) 162
\semi V. Pasquier, {\it J. Phys.} {\bf A20} (1987) 5707.}

\lref\VPae{V. Pasquier, {\OE nology of IRF models}, 
{\it Comm. Math. Phys.} {\bf 118} (1988) 355. }

\lref\Ko{I.K. Kostov, {\it Nucl. Phys.} {\bf B 300} [FS22] (1988) 559.}

\lref\DFZun{P. Di Francesco and J.-B. Zuber,
{\it Nucl. Phys.} {\bf B338} (1990) 602
.}
\lref\DFZ{P. Di Francesco and J.-B. Zuber,
in {\it Recent Developments in Conformal Field Theories}, Trieste Conference
1989, S. Randjbar-Daemi, E. Sezgin and J.-B. Zuber eds., World Scientific
1990 \semi
P. Di Francesco, {\it Int. J. Mod. Phys.} {\bf A7} (1992) 407.}

\lref\So{N. Sochen, {\it Nucl. Phys.} {\bf B360} (1991) 613. }

\lref\Ca{J. Cardy, {\it Nucl. Phys.} {\bf B324 } (1989) 581.}
\lref\SaBa{H. Saleur and M. Bauer, {\it Nucl. Phys. } {\bf B320} (1989) 591.}
\lref\PaSa{V. Pasquier and H. Saleur, {\it Nucl. Phys.} {\bf B330} (1990) 
523.}
\lref\PZ{V.B. Petkova and J.-B. Zuber, 
{\it Nucl. Phys.} {\bf B438} (1995) 347. }
\lref\PZdeu{ V.B. Petkova and J.-B. Zuber, 
\hepth 9510175, {\it Nucl. Phys. B} {\bf B463} (1996) 161. }

\lref\Verl{E. Verlinde, {\it  Nucl. Phys.} {\bf B300} [FS22] (1988) 360. }

\lref\We{H. Wenzl,
{\it Representations of Hecke Algebras and Subfactors},
 Ph.D. Thesis, University of Pennsylvania.}

\lref\Ocn{A. Ocneanu, {\it Quantized groups, string algebras
and Galois theory for algebras} in {\it Operator Algebras and Applications}, 
D.E. Evans and M. Takesaki edrs, 
Lond. Math. Soc. Lecture Note Series {\bf 136}, Cambridge University Press,
 119.}

\lref\Roc{P. Roche, {\it Comm. Math. Phys.} {\bf 127} (1990) 395.}

\lref\PZh{P. Pearce and Y.-K. Zhou, {\it Int. J. Mod. Phys. B} {\bf 7}
(1993) 3649.}

\lref\IZu{C. Itzykson and J.-B. Zuber, {\it Int. J. Mod. Phys.}
{\bf A7} (1992) 5661 and earlier references therein.}

\lref\SL{S. Loesch, {\it The exact solution of quasi one-dimensional
Hecke algebra models}, Diplomarbeit in Physik, Universit\"at Bonn report, 
November 1995.}

\Title{\vbox{\hbox{MRR 030-96}\hbox{SPhT 96/091 }
\hbox{{\tt q-alg/9608004  }}}}
{{\vbox {
\vskip-10mm
\centerline{ Identities in Representations of the Hecke Algebra    }
}}}
\medskip
\centerline{S. Loesch}\medskip
\centerline{ \it C.E.A.-SACLAY, Service de Physique Th\'eorique,} 
\centerline{ \it F-91191 Gif sur Yvette Cedex, France}
\bigskip
\centerline{Y.-K. Zhou}\medskip
\centerline{\it School of Mathematical Sciences, ANU,}
\centerline{\it Canberra ACT 0200, Australia}
\bigskip
\centerline{and}
\medskip
\centerline{J.-B. Zuber}\medskip
\centerline{ \it C.E.A.-SACLAY, Service de Physique Th\'eorique
,} \centerline{ \it F-91191 Gif sur Yvette Cedex, France}

\vskip .2in

\noindent 
New identities on traces of representations of the 
Hecke algebra on the spaces of paths on graphs are
presented. These identities are relevant in the 
computation of partition functions with fixed boundary
conditions and of two-point correlation functions in 
a class of lattice models. They might also
help to determine the Boltzmann weights of these (allegedly)
integrable models.

\bigskip

\Date{6/96\qquad\qquad submitted to  Int. J. Mod. Phys. }
%


\newsec{Introduction }
\noindent
This paper is devoted to the study of identities satisfied by the
Boltzmann weights of some two-dimensional integrable lattice models. 
These models are originally built from certain representations of 
the $sl(k)$ algebras, or rather from their quantum deformations. 
Imposing that their (`face') transfer matrices commute with the action 
of the quantum group ${\cal U}_qsl(k) $ implies that these
transfer matrices are generated by some representation of 
the Hecke algebra (in fact a quotient of the Hecke algebra by a condition 
depending on $k$, see below) \Res. 
Thus the original problem of building integrable models of that type is 
rephrased as the algebraic study of representations of the Hecke algebra.

In the graph models considered here, that are known  to have 
a critical regime described by minimal $W_{k}$ conformal theories, 
these representations are supported by the space of paths 
on a graph \VPun.  
Which graph is adequate for this construction remains to be fully
understood, but one may  start from the well-known `basic' case 
of the weight lattice of $sl(k)$ truncated at some level. 
A few other graphs have been discovered in previous investigations
\refs{\Ko\DFZun\DFZ{--}\So}, 
more have been conjectured lately \PZdeu, but the determination of their
Boltzmann weights and the verification that they satisfy  the relations of the
Hecke algebra is a painful, case-by-case analysis \So. 
It has been observed recently, however, that these weights seem to satisfy
general identities: traces of products of these weights take 
a simple, universal (i.e. graph independent)
expression. The aim of the present work was thus to understand 
the origin and the form of these identities and possibly to exploit them
for a general determination of the Boltzmann weights. 

The present work is also connected with considerations 
on conformal field theories. On the one hand, the 
traces that we are considering yield particular contributions 
to the partition functions of the lattice models 
with fixed boundary conditions. Along the line explored by
previous authors \refs{\SaBa,\DFZ}, these partition functions
are known to enjoy algebraic properties that match
those of the corresponding conformal field theories \Ca. 
On the other hand, these traces may also be regarded as 
particular contributions to two-point functions of certain 
lattice operators. The universality properties of correlation functions
with respect to the choice of representation of the Hecke algebra
have also consequences on properties of structure constants 
of the operator algebra of the corresponding conformal field theories: 
see \PZ\ for a discussion in the case of $sl(2)$. 

In sect. 2 we recall the context and notations following 
previous references. We also state the main results and 
conjectures of this paper. Section 3 presents the elements of
proof that we have collected.

%
\newsec{The Hecke Algebra. Identities on traces. }
\subsec{Notations on the Hecke algebra} 
\nind  
The Hecke algebra $\CH_{L}(q)$ is a $q$-deformation of the algebra of the 
symmetric group $S_L$. It is generated by the identity $1$ and by
 $g_i$, $i=1,\cdots, L-1$, 
the $q$-analogues of the transpositions $(i,i+1)$ generating 
the symmetric group $S_{L}$, that satisfy
\eqn\IIa
{\eqalign{g_i^2 &= 1 +(q-\qm) g_i\cr
g_i g_j& = g_j g_i \qquad |i-j|\ge 2; \cr
g_i g_{i+1} g_i &=  g_{i+1} g_i g_{i+1} \ .  \cr }}
Alternatively $\CH_L(q)$ is generated by the  generators $U_i=q-g_i$
satisfying
\eqn\IIb
{\eqalign{U_i^2 &= (q+\qm) U_i\cr
U_i U_j& = U_j U_i \qquad |i-j|\ge 2; \cr
U_i U_{i+1} U_i-U_i &=  U_{i+1} U_i U_{i+1} -U_{i+1}  \cr }}

We are interested in  representations of a quotient of the 
Hecke algebra by an additional condition. This comes about in 
the search of lattice integrable models whose diagonal-to-diagonal
transfer matrix commutes with the action of the quantum group
$\CU_qsl(k)$. The commutant of that quantum group in the
representation  space $[k]^{\otimes L}$ of the $L$-th power of the
$k$-dimensional fundamental representation is generated by
a certain $k$-dependent quotient of the Hecke algebra $\CH_{L}(q)$.

The extra condition relevant to describe the commutant 
of $\CU_q sl(k)$ is easier to express in terms of the $g$ generators. 
For any  permutation $\sigma$  of $S_{k+1}$, write a 
minimal decomposition into transpositions of adjacent elements:
$\sigma=\prod_{i\in I_{\sigma}} \sigma_i$ and denote by $g_{\sigma}$
the corresponding product of generators of the Hecke algebra
$g_{\sigma}=\prod_{i\in I_{\sigma}} g_i$. Then the condition 
reads \Res
\eqn\IIc{\sum_{\sigma\in S_{k+1}} (-q)^{-|I_{\sigma}|}  g_{\sigma}=0\ .}
See also \BaRe\ for a discussion of this condition.

\subsec{Representations on path spaces of graphs}
\nind We consider representations of this (quotient of the) Hecke
algebra on the space of paths of graphs. 

When dealing with $\CU_q sl(k)$ 
the basic example of graph 
is given by the system of integrable weights of the affine algebra
 $\widehat{sl(k)}_{n-k}$ at a certain level $n-k$, if  
$q=\exp {i\pi\over n}$. This system of weights is the
truncated weight lattice shifted by the vector $\bun:=\sum_1^{k-1}\GL_\l$
\eqn\IIca
{\CP_{++}^{(n)} =\{\Gl=\Gl_1\GL_1+\cdots 
+\Gl_{k-1}\GL_{k-1}|\Gl_\Ga \in \IN\,, \
\Gl_\Ga\ge 1,\ \Gl_1+\cdots \Gl_{k-1}\le n-1\} \ .}
We also need  the $k$ (linearly dependent) vectors $e_i$
\eqn\IIcb{\eqalign{
e_1=\GL_1,\qquad e_\Ga=\GL_\Ga-\GL_{\Ga-1},\quad i=2,\cdots
,&\ k-1,\qquad e_k=-\GL_{k-1} \cr
e_1+e_2+\cdots +e_k & =0\ , \cr }}
endowed with the scalar product
\eqn\IIcc
{(e_\Ga,e_\Gb)= \Gd_{\Ga\Gb}-{1\over k}\qquad\qquad \Ga,\Gb=1,\dots k\ .}
We then consider the fusion coefficients $N_{\Gl\Gm}^{\ \ \Gn}$ 
\Verl\ and the associative, 
commutative algebra that has these coefficients as structure constants.
The set of matrices $N_{\Gl}$, of non negative integer entries 
$(N_\Gl)_\Gm^{\ \Gn}=N_{\Gl\Gm}^{\ \ \Gn}$, for $\Gl,\Gm,\Gn \in 
\CP_{++}^{(n)}$, forms the regular representation of this fusion algebra
\eqn\IIcca{N_\Gl N_\Gm=N_{\Gl\Gm}^{\ \ \Gn}N_\Gn\ .}
If we denote for short 
$G_\l=N_{\bun+\GL_\l}$ the  fusion matrix for 
the $\l$-th fundamental representation of weight $\GL_\l$
(or rather $\bun+\GL_\l$ in the notation of \IIca),
it is known  that the set 
$\{N_\Gl\}$ is polynomially generated by the $G$'s
\eqn\IIccb{N_\Gl=\CP_\Gl(G_1,G_2,\cdots,G_{k-1})\ .}

The matrices $G_\l$, $\l=1, \cdots, {k-1}$ may be regarded as the adjacency 
matrices of a collection of 
graphs $\CG_\l$. Only $\CG_1$ will be relevant for the 
considerations of the present paper.  
In that basic case we shall refer to it as $\CA^{(n)}$: 
it is defined by its vertices, the weights of \IIca, and its 
edges between neighbouring vertices on the weight lattice oriented 
along the vectors $e_i$. 

\medskip
Beside this ``basic" case, the other graphs $\CG_\l$ 
that are known may be recast in that same framework by 
noticing that their matrices $G_1,G_2,\cdots,G_{k-1}$
generate by the same polynomials as in \IIccb\ another integrally valued
representation $ N^{(\CG)}_\Gl$ of the fusion algebra \IIcca\ (see \DFZun, 
where these matrices $\(N^{(\CG)}_\Gl\)_{ab}$ were denoted $(V_\Gl)_{ab}$.)

In the following, we shall also extend the notation $G_\l$ to
$G_0=G_k= \II$, (the identity matrix),  $G_\l=0$ if $\l>k$. 

\medskip
One then considers the space $\CP^{(L)}$ of paths of length $L$
on the graph $\CG_1$,
i.e. the space spanned by the states $\{|a_0 \cdots a_L \ket\}$
\eqn\IIcd{|a_0 \cdots a_L\ket
=(G_1)_{a_0 a_1}\,(G_1)_{a_1 a_2}  \cdots (G_1)_{a_{L-1} a_L}\,
|a_0\ket\otimes \, |a_1\ket  \otimes\, \cdots |a_L\ket\,,}
where $\bra a|b\ket=\Gd_{a b}\,$ and  the matrix elements of the
adjacency matrix $G_1$ of the graph $\CG_1$
ensure that consecutive vertices $a_i,a_{i+1}$ along the path are adjacent
on the graph.  (In this paper, for the sake of simplicity, 
we consider only graphs $\CG_1$ with no double edge: $(G_1)_{ab}=0$
or $1$.) 

When the graph under consideration is one of the basic ones
introduced above, it has been shown that the space of 
paths $\CP^{(L)}$ supports a representation
of the quotient of the Hecke algebra $\CH_L(q)$
given by the following formula \We
\eqnn\IIce
$$\eqalignno{ U_i| a_0 a_1 \cdots a_i \cdots a_L\ket
&=\sum_{a'_i}\,\, 
| a_0 a_1 \cdots a'_i \cdots a_L\ket \ 
\bdiamond{a_{i-1}}{a'_i}{a_{i+1}}{a_i} \ ,\qquad  i=1,\cdots, L-1\, ,
\cr
{\rm with \qquad }\bdiamond{a_{i-1}}{a_i}{a_{i+1}}{a'_i}
&= (1-\Gd_{\Ga\Gb}) \ 
{\[ s_{\Ga\Gb}(a_i) s_{\Ga\Gb}(a'_i)\]^{\oh}\over s_{\Ga\Gb}(a_{i-1})}
.&\IIce \cr}$$
where $a_i=a_{i-1}+e_\Ga$, $a_{i+1}=a_i+e_\Gb$ and the function $s_{\Ga\Gb}$
is defined as $s_{\Ga\Gb}(b)=\sin[(\pi/n)(e_\Ga-e_\Gb,b)]$. 
Here and in the following, the operators $U$ are regarded
as acting on vertical chains $a_0 \cdots a_L$
and {\it from right to left}, in agreement with the usual convention 
on composition of operators.

For other graphs, subject to additional conditions that are not
yet fully clarified, the space of paths also supports a representation of
the Hecke algebra. Only in the case of $\widehat{sl(2)}$ in which 
the relevant graphs turn out to be $ADE$ Dynkin diagrams, 
the general expression of the $U$'s has been found in terms of
the Perron-Frobenius eigenvector $\psi^{(\bun)}$ of $G_1$ \VPun
\eqn\IIcf{\bdiamond{a}{b}{c}{b'}= \Gd_{ac} {\(\psi^{(\bun)}_b
\psi^{(\bun)}_{b'}\)^{\oh} \over \psi^{(\bun)}_a} \ .}


\subsec{Traces}
\nind
Given a representation of the Hecke algebra of the previous type, 
namely on paths on a graph $\CG_1$, 
we need the trace of an arbitrary  element $x$ of that algebra. As the 
elements of the algebra do not affect the end points of the path, 
we may in fact compute the traces in the subspaces of paths 
with fixed ends $a_0$ and $a_L$. 
The trace $\tr$ of the element $x$ of the Hecke algebra $\CH_{L}$
with fixed $L$ is thus a matrix defined as 
\foot{The normalization constant
is sometimes fixed to be 
$1/(\Gc^{(1)}_1)^L$ so as to make the ``Markov trace"
$\tr_M(x):= \sum_{a_0,a_L}\psi^{(\bun)\,*}_{a_0}\psi^{(\bun)}_{a_L}
\(\tr(x)\)_{a_0\,a_L} $ independent of $L>L'$, if $x$ belongs to the 
algebra generated by $U_1,\cdots, U_{L'-1}$. Here and in 
the following, $\psi^{(\bun)}$ stands for the Perron-Frobenius eigenvector
common to all the $G$ matrices. We shall not make this
choice in this paper.}
\eqn\IIe{\(\tr(x)\)_{a_0\,a_L} := \sum_{a_1,a_2,\cdots,a_{L-1}} 
\langle a_0 a_1 \cdots a_{L-1}a_L|x| a_0 a_1 \cdots a_{L-1}a_L\rangle \ .}
In words, we are computing the trace of the element $x$ in the subspace 
$\CP^{(L)}_{a_0a_L}$ of paths of length
$L$ going from height $a_0$ to height $a_L$. 
The most important property of the trace \IIe\ that we shall need below is 
its cyclicity~: $\tr(x\,y)=\tr(y\,x)$. 

We are mainly interested in computing expressions of the form
\eqn\IId{\eqalign{(Z_L)& =\tr  U_1 U_2 \cdots  U_{L-1} 
\cr
{\rm or} \qquad (\Zt_L)& =\tr g_1 g_2 \cdots g_{L-1} 
\ , \cr
}}
which may be represented pictorially as in Fig. 1, with each 
diamond standing for a $U$ or a $g$, depending on the case. 
%
%
\fig{Pictorial representation of \IId. 
A summation over $a_1,\cdots, a_{L-1}$ is understood}
{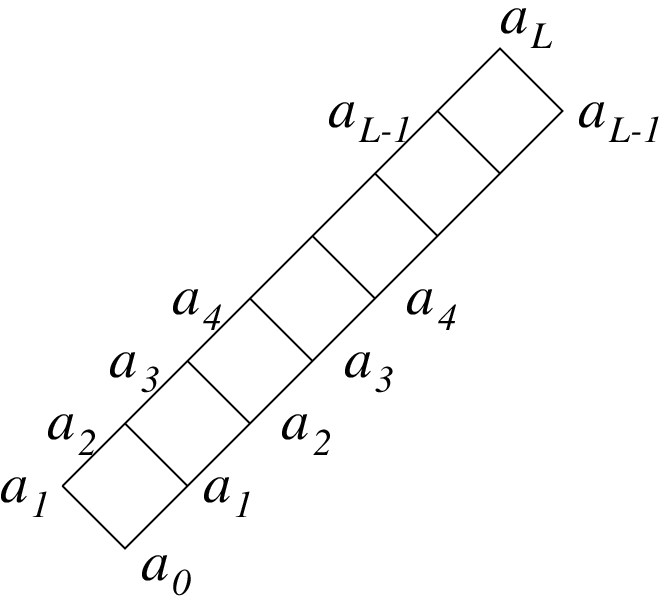}{4cm}\figlabel\squares
We recall that the local face transfer matrix at site $i$ of the integrable 
model built out of this representation of the Hecke algebra is a linear
combination of $U_i$ (or $g_i$) and the identity operator \VPun.  Thus
expressions  \IId\ are particular cases of what has to be 
computed in the evaluation of the partition function of the model. 

If $L=1$, the element $x$ is necessarily the identity, the summation
in \IIe\ is void and we set $Z_1=\Zt_1=G_1$ in agreement with \IIcd. 
Also it is easily verified that the computation of traces 
where some of the $U$'s or $g$'s are skipped reduces to the matrix 
product of  lower ones times possibly $G$'s; for example
\eqn\IIdaa{\eqalign {\(\tr U_1 U_3\)_{a_0 a_4} 
&= \sum_{a_1,\cdots a_3,a'_1,\cdots a'_3} 
\bra a_0 \cdots a_4|U_1|a_0 a'_1\cdots a'_3 a_4\ket 
\bra a_0 a'_1\cdots  a'_3 a_4 | U_3|a_0 \cdots a_4 \ket \cr
&= \sum_{a_1,a_2, a_3} 
\bra a_0  a_1 a_2 | U_1|a_0 a_1 a_2 \ket  
\bra a_2 a_3 a_4|U_3|a_2 a_3 a_4\ket 
=(Z_2. Z_2)_{a_0 a_4}\cr
}}
since $U_i$ ``sees" only heights $a_{i-1}$, $a_i$ and $a_{i+1}$ and 
affects only $a_i$. Likewise $\tr U_1 U_4= Z_2 G_1 Z_2$. A more complicated 
example is given by (Fig.2). 
\eqn\IIda{ \tr (U_2 U_4 U_5 U_8 U_{11}) = Z_1 Z_2 Z_3 Z_1 Z_2 Z_1^2 Z_2 
 \ .}
Notice that the labelling by $L$ of $Z_L$, $\Zt_L$ implies that
we are dealing with paths of length $L$ from $a_0$ to $a_L$, and 
the expression \IIda\ has indeed the grading $13$, if each $Z_\l$
is counted for $\l$. 
%
%
\fig{Pictorial representation of \IIda. 
The squares in solid line indicate the action of a $U$, those in doted 
line stand for the action of $1$. }{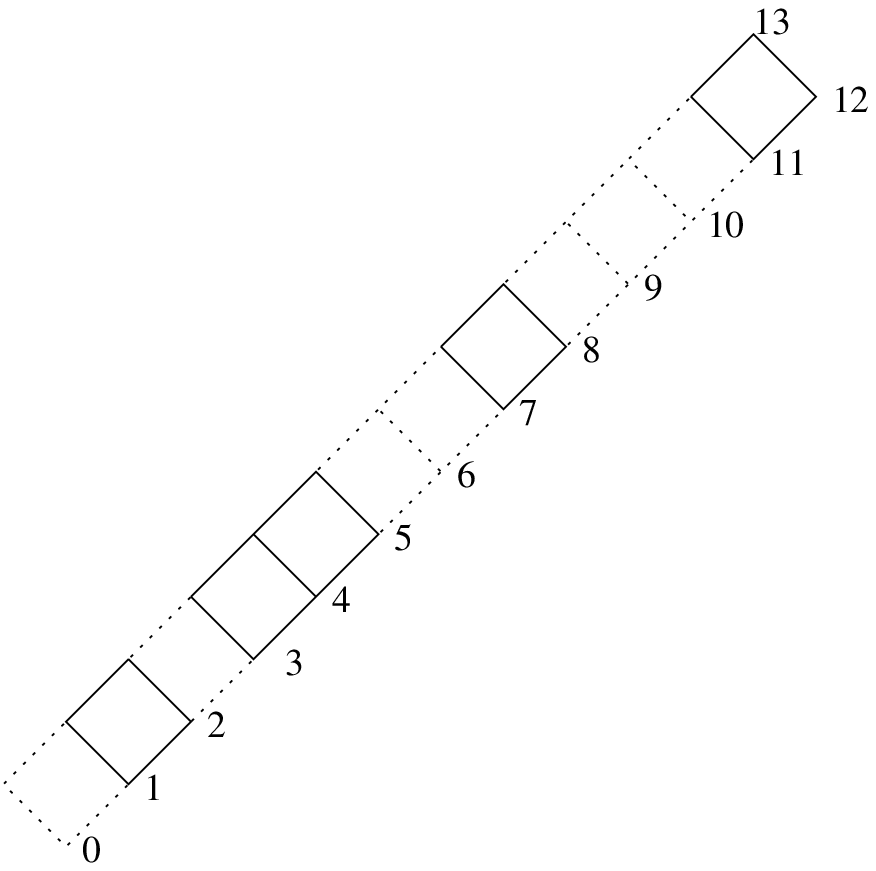}{4cm}\figlabel\mixed
%

\subsec{Properties of the $Z$ and $\Zt$ }
\nind

We show now that $Z$ and $\Zt$ are matrices that commute with 
the matrices $G_\l$ and may be expressed in terms of them. Moreover, we 
establish for low values of $k$ and $L$,  and conjecture in general, 
some recursion formulae  on $L$ for these expressions.

\medskip\noindent
{\bf a)} First we claim that the $Z$ and $\Zt$ matrices satisfy
\eqn\IIdd{Z^{(\CG)}_L = \sum_{\Gl\in \CP_{++}^{(n)}} z_L^{(\Gl)}
N^{(\CG)}_{\Gl} }
with the same matrices $N_\Gl^{(\CG)}$ as introduced above and 
a {\it universal} coefficient $z_L^{(\Gl)}$, and a similar 
relation involving coefficients  ${\tilde z}_L^{(\Gl)}$ for $\Zt_L$.
Here as above, `universal' means independent of the graph $\CG_1$ and
depending only on $k,n,L$; we shall see below that the dependence 
of $z_L^{(\Gl)}$ on $k,n$ is just embodied in the restriction of 
the summation over $\Gl$.

Consider first the case of the ``basic" graphs $\CA^{(n)}$. Then 
the matrix $N^{(\CA)}_\Gl$ is nothing else than the fusion 
matrix $N_\Gl$. As the proof of \IIdd\ is a bit lengthy and 
technical, it has been relegated to an appendix.


As mentionned in the Introduction, 
the fact that the traces \IId\ may be reduced through fusion 
matches the same property of partition functions with fixed 
boundary conditions in conformal field theories \Ca. 
Note that relation \IIdd\ allows to determine the coefficient $z_L^{(\Gl)}$
as
\eqn\IIdw{z_L^{(\Gl)} = (Z_L^{(\CA)})_{\bun \, \Gl} \ ,}
for the basic graph, and a similar relation for $\Zt_L^{(\Gl)}$. 
This is because $N_{\Gl\bun }^{\ \ \Gm}=\Gd_{\Gl\Gm}$. 
Thus \IIdd\ and \IIdw\ imply
\eqn\IIdx{ (Z_L^{(\CA)})_{\Gm \, \Gn} = (Z_L^{(\CA)})
_{\bun \Gl} 
N_{\Gl\Gm}^{\ \ \Gn} \ ,} 
The algebraic interpretation of \IIdx\ is that the only irreducible
representations of the Hecke algebra are those on paths on the basic graph 
running from $\bun$ to some weight $\Gl$ and that the fusion coefficients 
$N_{\Gm\Gn}^{\ \ \Gl} $ yield the multiplicities in the decomposition 
of the more general representation on paths from $\Gm $ to $\Gn$
 \refs{\PaSa,\Ca}.

For other graphs $\CG_1$, we have no direct proof of the general 
validity of \IIdd, but convergent evidences. First it can be 
proved for $k=2$, using the explicit Boltzmann weights \IIcf.  
Also as we shall see below, for $k=3$, very explicit forms 
will be given to the $Z_L$ based on the only defining relations of
the Hecke algebra, thus valid for all graphs, that agree with \IIdd. 
Finally the same interpretation may be given to the integer entries
$\(N_\Gl^{(\CG)}\)_{ab}$ of the matrices $N_\Gl^{(\CG)}$ 
as multiplicities in the decomposition 
of the representation carried by the paths on $\CG_1$ running between 
$a$ and $b$
\eqn\IIdx{ \(Z^{(\CG)}_L\)_{a \, b} = 
\(Z^{(\CA)}_L\)_{\bun \Gl} \(N^{(\CG)}_{\Gl}\)_{ab} \ .} 
%
%
\medskip\noindent
{\bf b)} Eq. \IIccb\ and \IIdd\ 
imply that $Z_L$  and $\Zt_L$ 
may be expressed in terms of the fundamental matrices 
$G_\l$, $\l=1,\dots,k-1$. In fact we have formulae
that enable us to determine these expressions recursively. 

For this we make use of Chebychev polynomials
of second type
\eqn\IIf{P_\l(\Gb) = {\sin (\l+1)\theta \over \sin \theta}\qquad 
\Gb=2\cos\theta.}
(We refrain from using the ordinary notation $U_\l$
as it could cause confusion with the generators of \IIa). 
They start as $P_0=1$, $P_1=\Gb$ etc. The recursion formula they 
satisfy, namely $\Gb P_\l=P_{\l+1}+P_{\l-1}$, may be used to extend 
them to $P_{-1}=0$, $P_{-2}=-1$. We consider 
these polynomials for $\theta$ given by $q=e^{i\theta}$, hence
$\Gb=q+\qm$. 

With these conventions we claim that
$\forall L>0  $
\eqna\IIg
$$ \eqalignno
{\sum_{\l=0}^{\inf(k,L-1)} (-q)^\l \Zt_{L-\l} G_\l & = (-1)^{L-1} 
P_{L-1}(\Gb) G_L  &\IIg a\cr
 \sum_{\l=0}^{\inf(k,L-1)} P_{\l-2}(\beta) Z_{L-\l} G_\l & =
-P_{L-1}(\Gb) G_L \ .   &\IIg b\cr}$$
Similar recursion formulae for the physical transfer matrix 
have been obtained in \SL. 
Eq. \IIg{a} or \IIg{b} define recursively all $\Zt_L$, resp $Z_L$.  
With our convention that $G_\l=0$ for $\l>k$, the right hand side 
of both eqs \IIg{} vanishes for $L>k$. 
Thus \IIg{} reads more explicitly
\eqna\IIga
$$  \eqalignno
{L \ge k+1 \qquad \Zt_L &= - \sum_{\l=1}^{k} (-q)^{\l} \Zt_{L-\l} G_\l 
 &\IIga a\cr
L\le k \qquad \Zt_L &= - \sum_{\l=1}^{L-1} (-q)^{\l} \Zt_{L-\l} G_\l +
(-1)^{L-1}  P_{L-1}(\Gb) G_L  &\IIga b\cr}$$
and likewise
\eqna\IIgb
$$ \eqalignno
{L\ge k+1 \qquad Z_L &= \sum_{\l=2}^{k}P_{\l-2}(\beta) Z_{L-\l} G_\l 
 &\IIgb a\cr
L\le k \qquad  Z_L &= \sum_{\l=2}^{L-1} P_{\l-2}(\beta) Z_{L-\l} G_\l +
P_{L-1}(\Gb) G_L \ .  &\IIgb b\cr}$$

Note that because of the vanishing of $P_1$, the dependence 
of $Z_L$ on $G_1$ comes only through that of $Z_1,\cdots,Z_k$, 
and $Z_L$ is thus at most of degree one in $G_1$.

The first $Z$ and $\Zt$ read
\eqn\IIh{\eqalign{
Z_1 &= G_1 \cr
Z_2 &= \beta G_2 \cr
Z_3 &= G_1 G_2 + (\Gb^2-1) G_3 \cr
Z_4 &= \Gb\(G_2^2+G_1G_3\) + \Gb(\Gb^2-2) G_4 \cr
Z_5 &= G_1 G_2^2 + (2\Gb^2-1) G_2G_3 +(\Gb^2-1)G_1G_4 +(\Gb^4-3\Gb^2+1) G_5\cr
 \vdots &\qquad \qquad\qquad\vdots \cr}}
\eqn\IIi{\eqalign{
\Zt_1 &= G_1 \cr
\Zt_2 &= q G_1^2 - (q+\qm) G_2 \cr
\Zt_3 &= q^2 G_1^3 -(2q^2+1) G_1 G_2 + (q^2+1+q^{-2}) G_3 \cr
\Zt_4 &= q^3 G_1^4 -(3q^3+q) G_1^2 G_2 + (q^3+q) G_2^2 
+ (2q^3+q+\qm) G_1G_3 
\cr & \qquad - (q^3+q+\qm+q^{-3}) G_4 \cr
 \vdots &\qquad \qquad\qquad \vdots \cr}}
and the reader recognizes in the coefficient of $G_L$ in $Z_L$
resp. $\Zt_L$ two alternative expressions of $P_{L-1}(\Gb)$ (up
to a sign).
%
%
%
%
\def\tvp{\vrule height 2pt depth 1pt} 
\def\thp{\vrule height 0.4pt width 0.35em}
\def\cc#1{\hfill#1\hfill}

\setbox201=\vbox{\offinterlineskip
\cleartabs
\+ \thp&\cr  
\+ \tvp\cc{}&\tvp\cr 
\+ \thp&\cr          
\+ \tvp\cc{}&\tvp\cr 
\+ \thp&\cr  }

\setbox3=\hbox{$\vcenter{\offinterlineskip
\+ \thp&\cr  
\+ \tvp\cc{}&\tvp\cr 
\+ \thp&\cr  	     
\+ \tvp\cc{}&\tvp\cr 
\+ \thp&\cr  
\+ \tvp\cc{}&\tvp\cr 
\+ \thp&\cr  }$}

\setbox4=\vbox{\offinterlineskip
\cleartabs
\+ \thp&\cr  
\+ \tvp\cc{}&\tvp\cr 
\+ \thp&\cr  	     
\+ \tvp\cc{}&\tvp\cr 
\+ \thp&\cr  	     
\+ \tvp\cc{}&\tvp\cr 
\+ \thp&\cr  
\+ \tvp\cc{}&\tvp\cr 
\+ \thp&\cr  }

\setbox5=\vbox{\offinterlineskip
\cleartabs
\+ \thp&\cr  
\+ \tvp\cc{}&\tvp\cr 
\+ \thp&\cr  	     
\+ \tvp\cc{}&\tvp\cr 
\+ \thp&\cr   	     
\+ \tvp\cc{}&\tvp\cr 
\+ \thp&\cr  
\+ \tvp\cc{}&\tvp\cr 
\+ \thp&\cr  
\+ \tvp\cc{}&\tvp\cr 
\+ \thp&\cr  }

\setbox21=\vbox{\offinterlineskip
\cleartabs
\+ \thp&\thp&\cr  
\+ \tvp\cc{}&\tvp\cc{}&\tvp\cr 
\+ \thp&\thp&\cr  		
\+ \tvp\cc{}&\tvp\cr 
\+ \thp&\cr  }

\setbox31=\vbox{\offinterlineskip
\cleartabs
\+ \thp&\thp&\cr  
\+ \tvp\cc{}&\tvp\cc{}&\tvp\cr 
\+ \thp&\thp&\cr  
\+ \tvp\cc{}&\tvp\cr 
\+ \thp&\cr 	     
\+ \tvp\cc{}&\tvp\cr 
\+ \thp&\cr  }

\setbox113=\vbox{\offinterlineskip  
\cleartabs
\+ \thp&\thp&\thp&\cr  
\+ \tvp\cc{}&\tvp\cc{}&\tvp\cc{}&\tvp\cr 
\+ \thp&\thp&\thp&\cr  
\+ \tvp\cc{}&\tvp\cr 
\+ \thp&\cr 	     
\+ \tvp\cc{}&\tvp\cr 
\+ \thp&\cr  }

\setbox22=\vbox{\offinterlineskip
\cleartabs
\+ \thp&\thp&\cr  
\+ \tvp\cc{}&\tvp\cc{}&\tvp\cr 
\+ \thp&\thp&\cr  	       
\+ \tvp\cc{}&\tvp\cc{}&\tvp\cr 
\+ \thp&\thp&\cr}

\setbox41=\vbox{\offinterlineskip
\cleartabs
\+ \thp&\thp&\cr  
\+ \tvp\cc{}&\tvp\cc{}&\tvp\cr 
\+ \thp&\thp&\cr  	       
\+ \tvp\cc{}&\tvp\cr 
\+ \thp&\cr 		       
\+ \tvp\cc{}&\tvp\cr 
\+ \thp&\cr 
\+ \tvp\cc{}&\tvp\cr 
\+ \thp&\cr  }

\setbox32=\vbox{\offinterlineskip	
\cleartabs				
\+ \thp&\thp&\cr  			
\+ \tvp\cc{}&\tvp\cc{}&\tvp\cr 
\+ \thp&\thp&\cr  
\+ \tvp\cc{}&\tvp\cc{}&\tvp\cr 
\+ \thp&\thp&\cr
\+ \tvp\cc{}&\tvp\cr 
\+ \thp&\cr }

\setbox221=\vbox{\offinterlineskip
\cleartabs
\+ \thp&\thp&\thp&\cr  
\+ \tvp\cc{}&\tvp\cc{}&\tvp\cc{}&\tvp\cr 
\+ \thp&\thp&\thp&\cr  
\+ \tvp\cc{}&\tvp\cc{}&\tvp\cr 
\+ \thp&\thp&\cr}		

\setbox101=\vbox{\offinterlineskip		
\cleartabs				
\+ \thp&\cr  
\+ \tvp\cc{}&\tvp\cr 
\+ \thp&\cr  }

\setbox200=\vbox{\offinterlineskip	
\cleartabs
\+ \thp&\thp&\cr  
\+ \tvp\cc{}&\tvp\cc{}&\tvp\cr 
\+ \thp&\thp&\cr}
\setbox155=\vbox{\offinterlineskip
\cleartabs				
\+ \thp&\thp&\thp&\cr  			
\+ \tvp\cc{}&\tvp\cc{}&\tvp\cc{}&\tvp\cr 
\+ \thp&\thp&\thp&\cr    }
\setbox160=\vbox{\offinterlineskip
\cleartabs				
\+ \thp&\thp&\thp&\thp&\cr  			
\+ \tvp\cc{}&\tvp\cc{}&\tvp\cc{}&\tvp\cc{}&\tvp\cr 
\+ \thp&\thp&\thp&\thp&\cr    }
\setbox199=\vbox{\offinterlineskip
\cleartabs				
\+ \thp&\thp&\thp&\thp&\thp&\cr  			
\+ \tvp\cc{}&\tvp\cc{}&\tvp\cc{}&\tvp\cc{}&\tvp\cc{}&\tvp\cr 
\+ \thp&\thp&\thp&\thp&\thp&\cr   }
\setbox180=\vbox{\offinterlineskip
\cleartabs				
\+ \thp&\thp&\thp&\thp&\cr  		
\+ \tvp\cc{}&\tvp\cc{}&\tvp\cc{}&\tvp\cc{}&\tvp\cr 
\+ \thp&\thp&\thp&\thp&\cr   
\+ \tvp\cc{}&\tvp\cr           			   
\+ \thp&\cr 	    }
\setbox211=\vbox{\offinterlineskip  
\cleartabs
\+ \thp&\thp&\thp&\cr  
\+ \tvp\cc{}&\tvp\cc{}&\tvp\cc{}&\tvp\cr 
\+ \thp&\thp&\thp&\cr  
\+ \tvp\cc{}&\tvp\cr 
\+ \thp&\cr 	     
}
%

%
\medskip
\noindent
{\bf c)} We can also  rewrite the right hand sides of \IIh\ as 
{\it linear} expressions in the $N$'s, thus getting explicit 
expressions for \IIdd. Curiously, the first polynomials
in $\Gb$ turn out to have only positive coefficients. For example
the first ones read
\eqn\IIm{\eqalign{
Z_1 &= N_{\copy101} \cr
Z_2 &= \beta N_{\copy201} \cr
Z_3 &= \Gb^2  N_{\copy3} + N_{\copy21}  \cr
Z_4 &= \Gb^3 N_{\copy4}  + \Gb \( N_{\copy31} + N_{\box22} \)\cr
Z_5 &= \Gb^4 N_{\copy5} + \Gb^2 \(3 N_{\copy41} +2 N_{\copy32} \) +
 \(N_{\box32} + N_{\copy 113} +N_{\box221} \) \cr
 \vdots &\qquad \vdots \cr}}
We have checked that up to $Z_9$.
The reason of this intriguing feature is not known to us, but should 
have to do with positivity properties of the operator $U_1\cdots U_{L-1}$. 
It also suggests that the coefficients $z_L^{(\Gl)}$ might be 
given a combinatorial or  group theoretic interpretation.\foot{Here 
again, to get universal expressions, we  have
 extended slightly the original conventions on the 
weights or Young tableaux labelling the fusion matrices $N$: 
for a given $k$, any column in a  Young tableau with  $k$ rows may be 
deleted, any tableau with a column of more than $k$ rows does not
contribute. Thus for instance $N_{\copy4}=\II$ if $k=4$, and  $=0$ for 
$k\le 3$.}

Even simpler are the corresponding expressions for the $\Zt_L$
\eqn\IIma{\eqalign{
\Zt_1 &= N_{\copy101} \cr
\Zt_2 &= q  N_{\copy200}- q^{-1} N_{\copy201} \cr
\Zt_3 &= q^{2} N_{\copy155} - N_{\copy21} +q^{-2} N_{\copy3}   \cr
\Zt_4 &= q^{3} N_{\box160} -q N_{\box211} +q^{-1} N_{\box31} 
-q^{-3} N_{\box4} 
 \cr
\Zt_5 &=  q^{4} N_{\box199} -q^{2} N_{\box180} + N_{\box113} -q^{-2} 
N_{\box41} +q^{-4} N_{\box5} 
\cr }\ .}
These formulae suggest an obvious generalization. 
Consider the Young tableaux 
with at most one line of length larger than one 
(i.e. of type $(t+1) (1)^s$, $s+t+1=L$, $s,t \ge 0$) and denote by 
$\Gl_{t,s}$ the corresponding weight of $sl(k)$ ($k \ge s+1$)
\def\tv{\vrule height 4pt depth 2pt} 
\def\th{\vrule height 0.4pt width 0.7em} 
\def\ca{\hfill{\ }\hfill}
\setbox111=\hbox{$\vbox{\offinterlineskip
\+          &&   &\th&\th&\th&\th&\th&\th&\th&\th&\th&\cr
\+ 
&\tv&\ca&\tv\ca&\tv\ca&\tv\ca&\tv\ca&\tv\ca&\tv\ca&\tv\ca&\tv\ca&\tv\ca&
\tv&\cr
\+       &  &&\th&\th&\th&\th&\th&\th&\th&\th&\th&\cr
}$}
\setbox12=\hbox{$\vcenter{\offinterlineskip
\+       &\tv&\ca&\ca&\tv&   &   &   &   &   &\cr
\hrule
\+       &\tv&\ca&\ca&\tv&   &   &   &   &   &\cr
\hrule
\+ \cr
\+        & \tv&\ca&\ca&\tv&   &   &   &   &   & \cr
\hrule }$}
\setbox22=\hbox{$\left.\vbox to \ht12{}\right\}$}
\setbox33=\hbox{$s$}
\setbox10=\hbox to 6.2em {\downbracefill}
$$\Gl=\Gl_{t,s}:=\quad \vcenter{\vbox{\offinterlineskip\halign{
#& \qquad #& \hfill # \cr
\hbox to 6.5em{\hfill\quad$ t+1$ \quad\hfill}\cr
\noalign{\vskip 2mm}
\box10 \cr
\noalign{\vskip 2mm}
\copy111 \cr
\copy12\copy22\box33 \cr  }} }$$
Then 
\eqn\IImb{\Zt_L=\sum_{{\Gl_{t,s}\atop  s+t+1=L} }
  (-1)^s q^{L-1-2s} N_{\Gl_{t,s}} \ .}
As we shall see in sect. 3.1, this conjectured expression
 is indeed the solution to 
eq. \IIg{a} or \IIga{} and it reduces to a well known fact 
in the classical limit $q\to 1$.

\bigskip
%
%
{\bf d)} Explicit solutions to eqs \IIg{b} may be given for $k=2,3$. 
For $k=2$, 
\eqn\IIn{\eqalign{Z_{2\l+1}&= G_1 G_2^\l \cr Z_{2\l}&= \Gb G_2^\l\ .\cr} }
For $k=3$,
\eqn\IIo{\eqalign{Z_{2\l+1}&=\sum_{p=0}^{[{\l\over 3}]}{\l-p\choose 2p}
\Gb^{2p}G_1 G_2^{\l-3p} G_3^{2p} \cr
& \qquad\qquad\qquad
+\sum_{p=0}^{[{\l-1\over 3}]}
\[{\l-p\choose 2p+1}\Gb^{2p+2}-{\l-p-1\choose 2p}\Gb^{2p}\]G_2^{\l-3p-1}
G_3^{2p+1} \cr
Z_{2\l}&= \Gb G_2^\l+ \sum_{p=1}^{[{\l\over 3}]}
\[{\l-p\choose 2p}\Gb^{2p+1}-{\l-p-1\choose 2p-1}\Gb^{2p-1}\]G_2^{\l-3p}
G_3^{2p} \cr
& \qquad\qquad\qquad
+\sum_{p=0}^{[{\l-2\over 3}]}{\l-p-1\choose 2p+1}
\Gb^{2p+1}G_1 G_2^{\l-3p-2} G_3^{2p+1}\cr 
} \ .}
In these two cases $k=2,3$, we have respectively $G_2\equiv 1$, 
$G_3\equiv 1$, but we have left them in order to make
 explicit the homogeneity of $Z_\l$. 
In checking that these expressions do satisfy the recursion 
formulae \IIgb{} it is useful to use the following formula
for $P_n(\Gb)$:
\eqn\IIp{P_n(\Gb)=\sum_{p=0}^{[{n\over 2}]}(-1)^p{n-p\choose p}\Gb^{n-2p}\ .}
Similar expressions may be obtained for higher $k$ but 
they become increasingly cumbersome. 
Eq. \IIn\ may be easily proved in full generality
for all graphs using Pasquier's expressions for the Boltzmann
weights \VPun.
\bigskip
%
%
{\bf e)} 
The status of the assertions \IIg{} is as follows
\item{1)} the expression \IImb\ may be proved to solve \IIg{a},  and its
limit as $q\to 1$ is related to well known properties of 
 the symmetric group; as a result, \IIg{a} itself is proved in that 
limit $q=1$  (see below, sect. 3.1);
\item{2)} eq. \IIg{a} may be shown to be a 
simple consequence of \IIc\ for $k=2,3$ and all $L>k$ (see below, sect. 3.2);
thus it holds for {\it all} graphs; 
\item{3)} eq. \IIg{b} has been checked for the basic graph 
of $sl(k)$ up to $k=9$ and $L=9$ using the explicit 
expression \IIce;
\item{4)} eq. \IIg{b} has been proved for the basic graph of 
$sl(k)$, $k=3,4$, all $L>0$ (see below, sect. 3.3);
\item{5)} 
the non trivial equivalence between \IIg{a} and \IIg{b} 
has been  verified for $k=2,3,4$;
\item{6)} these identities are consistent with another 
identity, related to the Markov property of the trace (see footnote 1), 
namely that 
\eqn\IIj{ \sum_c \bdiamond{a}{b}{c}{b}\ \psi^{(1)}_c = P_{k-2}(\Gb) 
\psi^{(1)}_b \  , }
which may be proved on the expression \IIce. (For other graphs, it 
then follows from the existence of `cells' \refs{\Ocn,\DFZun,\Roc,\PZh}). 
From that identity, it follows that
\eqn\IIk{\langle Z_L\rangle :=
\sum_{a,b} \psi^{(1)}_a (Z_L)_{ab} \psi^{(1)}_b=P_{k-1}(\Gb)
P^{L-1}_{k-2}(\Gb)\ ,}
and it is easy to see that \IIk\ satisfies the recursion 
formula \IIg{b}. 

In view of these evidences, it seems quite reasonable 
to conjecture the general validity of \IIg{} for all
representations, all $k$ and all $L>k$.

%
\newsec{Elements of proof. }
\noindent In this section we collect the various evidences 
on the validity of eqs \IIg{} that we have found
\subsec{\IImb\ solves \IIga{}; the $q\to 1$ limit of eq. \IIga{}} 
\noindent
First let us consider the expression \IImb\ and let us show why it 
is natural in the limit where $n\to \infty$, hence $q\to 1$. 
Then the paths we consider in the basic case are just those
on the entire weight lattice of $sl(k)$. 
For $q$, the deformation parameter equal to 1, the Hecke algebra
$\CH_L$ reduces to the algebra of the symmetric group $S_L$ 
and its irreducible representations associated with paths of length 
$L$ from $\bun$ to $\lambda$ are nothing else than the 
representations of $S_L$ attached to the Young tableau 
associated with $\Gl$. Thus for any $\Gs\in\CH_L$, the matrix element
$(\tr\,\Gs)_{\bun\Gl}$ is simply the irreducible character 
$\chi_\Gl(\Gs)$ of $S_L$. 
Moreover, the product $g_1 g_2 \dots g_{L-1}$
whose trace defines $\Zt_L$ in \IId\ is a cyclic permutation
of $S_L$. It is well known (see for example \IZu) that the only Young 
tableaux giving non vanishing characters of cyclic permutations $\Gs$ of $S_L$
are those denoted   $\Gl_{t,s}$ in sect. 2.4.c,  for  which
  $\chi_{\Gl_{t,s}}(\Gs)=(-1)^s$.
Thus in the limit $q=1$, 
\eqn\IIIaa{(\Zt_L)_{\bun \Gl} = \sum_{{\Gl_{t,s}} }
(-1)^s \Gd_{\Gl \Gl_{st}}}
or 
$$ {\Zt_L=\sum_{{\Gl_{t,s}\atop  s+t+1=L} }
  (-1)^s  N_{\Gl_{t,s}} \ ,}
    \eqno{\IImb_1} $$ 
and  expression \IImb\ is a very natural $q$-deformation of that formula. 

In fact one may prove that \IImb\ satisfies \IIga{} for all $q$. 
Considering first the case where $k \ge L$,  insert \IImb\ in 
the $(1,\Gm)$ matrix element of the l.h.s. of \IIg{a}. 
\eqn\IIIb{ \sum_\Gl
 \sum_{\l=0}^{L-1}(-q)^\l (\Zt_{L-\l})_{\bun \Gl} (G_\l)_{\Gl\Gm}
\buildrel{\rm ?}\over=   \sum_{\l=0}^{L-1} 
\sum_{{\Gl_{t,s}\atop  s+t+1=L-\l} }    (-1)^s  q^{L-1-\l-2s} (-q)^\l 
 (G_\l)_{\Gl_{t,s}\Gm} }
It receives contributions from those weights $\Gm$ that 
appear in the tensor product of $\Gl_{t,s}$ by the 
$\l$-th fundamental representation
\setbox5=\hbox{$\vcenter{\offinterlineskip
\hrule
\+       &\tv&\ca&\ca&\tv&   &   &   &   &   &\cr
\hrule
\+       &\tv&\ca&\ca&\tv&   &   &   &   &   &\cr
\hrule
\+       &\tv&\ca&\ca&\tv&   &   &   &   &   &\cr
\hrule
\+       &\tv&\ca&\ca&\tv&   &   &   &   &   &\cr
\hrule 
\+       &\tv&\ca&\ca&\tv&   &   &   &   &   & \cr
\hrule }$}
\setbox15=\hbox{$\left.\vbox to \ht5{}\right\}$}
\setbox25=\hbox{$\l$}
$$\vcenter{\offinterlineskip\halign{
#& \qquad #& \hfill # \cr
\box111 \cr
\copy12 \cr}}\otimes \vcenter{\hbox{\copy5\copy15\copy25}} $$
It is  then a simple exercise of bookkeeping to see that 
there are cancellations for almost all $\Gm$ and that only 
the vertical Young tableau with $L$ boxes contributes  
$(-1)^{L-1} (q^{L-1}+ q^{L-3}+ \cdots + q^{-(L-1)}) $ $= (-1)^{L-1} 
P_{L-1}(\Gb)$
 in agreement with the conjectured formula. 
For example, for $L=3$, $k\ge 3$
\def\tvp{\vrule height 2pt depth 1pt} 
\def\thp{\vrule height 0.4pt width 0.35em}
\def\cc#1{\hfill#1\hfill}
\setbox200=\vbox{\offinterlineskip
\cleartabs
\+ \thp&\thp&\cr  
\+ \tvp\cc{}&\tvp\cc{}&\tvp\cr 
\+ \thp&\thp&\cr}
\setbox100=\vbox{\offinterlineskip
\cleartabs
\+ \thp&\cr  
\+ \tvp\cc{}&\tvp\cr 
\+ \thp&\cr  }
\eqn\IIIc{
\eqalign{\Zt_3 N_0- q \Zt_2 N_1+q^2 \Zt_1 N_2&=
\(q^2 \copy155-\copy21+q^{-2} \copy3\) -q\(q\copy200-q^{-1}\copy201 \)\otimes 
\copy101+q^2 \copy101\otimes \copy201\cr
&=q^2\copy155-\copy21+q^{-2}\copy3
-q^2\copy155-q^2 \copy21+\copy21+\copy3
+q^2\copy21+q^2\copy3\cr
&=(q^2+1+q^{-2})\ \copy3\ .}}
For $k < L$, 
the same discussion applies, provided one restricts oneself to 
Young tableaux with less than $k+1$ rows. For instance,  
 for $k=2$, the contribution of $\copy3$ in \IIIc\ has to be removed.

In this section, we have proved that i) eq. \IImb\ is true for 
$q=1$; ii) that it satisfies \IIg{a} for all $q$. This implies that 
the conjectured identity \IIg{a} is itself established for $q=1$. 


%
\subsec{{\rm \IIc} implies  
\IIg{a}}
\noindent 
This implication has been checked for $k=2,3$. For $k=3$, to write explicitly 
\IIc, we have to enumerate all 
permutations of the symmetric group  $S_4$ and to decompose them into 
products of adjacent transpositions $g_i$, $i=1,2,3$. We get in that way
\eqn\IIId{\eqalign{ &
 I - \qm (g_1+g_2+g_3) +q^{-2} (g_1g_3+g_2g_1+g_1g_2+g_2g_3+g_3g_2)\cr
 &-q^{-3} (g_1g_2g_1+g_2g_3g_2+g_3g_2g_1+g_1g_3g_2+g_1g_2g_3+g_2g_3g_1) \cr
& +q^{-4}(g_2g_1g_3g_2+g_3g_1g_2g_3+g_3g_2g_1g_3+g_1g_2g_3g_1+g_1g_3g_2g_1)\cr
& -q^{-5}\,(g_2g_3g_2g_1g_2+g_2g_3g_1g_2g_3+g_1g_2g_3g_2g_1)
+q^{-6}\,  g_2g_1g_2g_3g_2g_1 =0 \ .}}
Multiply this by $g_4 g_5 \cdots g_{L-1}$, take the trace and 
use eqs \IIa\ and the cyclicity of the trace to simplify this expression. 
A straightforward but tedious calculation yields
\eqn\IIIe{ (1+q^{-2})(q^{-3}+q^{-5}+q^{-7})\(-\Zt_L+q \Zt_{L-1}G_1-q^2 
\Zt_{L-2}G_2 +q^3 \Zt_{L-3}G_3 \) =0 \ }
which is the appropriate eq. \IIg{a}. 
The same can be done for $k=2$ and might be repeated for all $k>3$, although
we have not been able to master the combinatorics.
%


%

%
%
\subsec{Proving \IIg{b} for all $L$ and $k=3,4$}
\noindent 
For the basic graph for $k=3,4$, using the explicit expressions \IIce\ of the
matrix elements of the $U$'s, we have been able to prove 
\IIg{b} for all $L\ge k$ (while it is verified for $L\le k$
by simple inspection). The proof displays an 
interesting and suggestive hierarchical structure. 

\def\arr#1{\buildrel {#1} \over {\mapsto}}
Let us sketch the proof for $k=3$. The identity to be proved 
reads 
\eqn\IIIf{ Z_L=Z_{L-2}+\beta Z_{L-3} \ .}
Denote $a \arr{\Ga} b$ 
an edge between the two heights $a$ and $b$ of the ``basic graph"
($a$ and $b$ are integrable weights)
along the direction $e_\Ga$. We look at the contributions 
to $Z_{L}$ from paths ending at some height $e$ that passed
through the heights $a$ and $b$, i.e. paths
\eqn\IIIg{ a \arr{\Ga} b \arr{\Gb} c \arr{\Gd} d  \arr{\Ge} e\ .}
with $a,b,\Ga,e$ fixed, and we shall sum over the possible
$c,d$. The direction $\Gb$ must be different from $\Ga$
and {\it by convention}, $\Gc$ denotes the third direction 
such that $e_\Ga+e_\Gb+e_\Gc=0$. 

Two cases are possible: 
\item{1)} $\Gb$, $\Gd$ and $\Ge$ are all different, hence equal
up to some permutation to $\Ga,\Gb,\Gc$, hence $e=b$. This is 
possible for
 $$ a \arr{\Ga} b \arr{\Gb} c \arr{\Ga}d  \arr{\Gc} e=b\ $$
or for
 $$ a \arr{\Ga} b \arr{\Gb} c \arr{\Gc}d  \arr{\Ga} e=b\ ;$$
\item{2)} $\Gb,\Gd,\Ge$ not all different, either as
 $$ a \arr{\Ga} b \arr{\Gb} c \arr{\Gc}d  \arr{\Gb} e=c+e_\Gb+e_\Gc\ $$
or as
 $$ a \arr{\Ga} b \arr{\Gb} c \arr{\Ga}d  \arr{\Gb} e=c+e_\Gb+e_\Ga\ .$$

Let us now compute the contributions of these paths to the l.h.s.
of the recursion formula
\eqn\IIIh{ \eqalign{
\sum_{c,d} Z^{(L-3)}_{ab} \bdiamond{a}{b}{c}{b}\ \bdiamond{b}{c}{d}{c}
\ \bdiamond{c}{d}{e}{d}&=\sum_j Z^{(L-3)}_{ab}
\bdiamond{a}{b}{\!\!\!\!\!\!b+e_\Gb}{b} \cr
\Big[  
\bdiamond{b}{c=b+e_\Gb}{\!\!\!\!\!\!\!\!\!d=c+e_\Gc}{\!c}\ \bdiamond{c}{d}{b}{d}
&+\bdiamond{b}{c=b+e_\Gb}{\!\!\!\!\!\!\!\!\!d=c+e_\Ga}{\!c}\ 
\bdiamond{c}{d}{b}{d}
 \cr
+\ \bdiamond{b}{c=b+e_\Gb}{\!\!\!\!\!\!\!\!\!d=c+e_\Gc}{\!c}\ 
\bdiamond{c}{d}{\!\!\!\!\!\!\!\!\!e=c+e_\Gb+e_\Gc}
{\!\!\!\!\!\!\!\!\!\!d}
&+\ \bdiamond{b}{c=b+e_\Gb}{\!\!\!\!\!\!\!\!\!d=c+e_\Ga}{\!c}
\ \bdiamond{c}{d}{\!\!\!\!\!\!\!\!\!e=c+e_\Gb+e_\Ga}{\!\!\!\!\!\!\!\!\!\!c}
\Big]   \cr } }
where $ Z^{(L-3)}_{ab} $ is  the contribution to
$Z_{L-2}$ from the paths ending along the edge $ab$. 

The last two terms in the square bracket  are equal to one. 
Indeed with the notations of \IIce
\eqn\IIIha{\bdiamond{b}{c=b+e_\Gb}{\!\!\!\!\!\!\!\!\!d=c+e_\Gc}{c}
\ \bdiamond{c}{d}{\!\!\!\!\!\!\!\!\!\!\!\!e=c+e_\Gb+e_\Gc}{d}
={s_{\Gb\Gc}(c)\over s_{\Gb\Gc}(b)} {s_{\Gc\Gb}(d)\over s_{\Gc\Gb}(c)}=1 }
since $c=b+e_\Gb$ hence $s_{\Gc\Gb}(d)=s\((e_\Gc-e_\Gb),(b+e_\Gc+e_\Gb)\)=-
s_{\Gb\Gc}(b)$ and $s_{\Gb\Gc}(c)=-s_{\Gc\Gb}(c)$.

Let us compare the result of \IIIh\ 
with the contributions to the r.h.s. of \IIIg\
of the same paths  through $a$ and $b$
\eqn\IIIi{ 
Z^{(L-2)}_{ab} \Big[
\sum_j \bdiamond{a}{b}{\!\!\!\!\!\!\!\!\!c=b+e_\Gb}{b}\ \ (G_2)_{ce}
+\beta \delta_{be}\Big] \ .}
Since $(G_2)_{ce}=\delta_{e,c+e_\Ga+e_\Gb}+\delta_{e,c+e_\Gb+e_\Gc}
+\delta_{e,c+e_\Gc+e_\Ga}$, the first two terms reproduce the last two in (2), 
and the last   is equal to $\delta_{eb}$, hence we
have to compare all the terms in \IIIh\ and \IIIi\ proportional to 
$\delta_{be}$. Let us show indeed that ($i$ being fixed, {\it not} summed 
over)
\eqn\IIIj{ \sum_{\Gb\ne \Ga} \bdiamond{a}{b}
{\!\!\!\!\!c=b+e_\Gb}{b=a+e_\Ga}
\Big[
\bdiamond{b}{c=b+e_\Gb}{\!\!\!\!\!\!\!\!\!d=c+e_\Gc}{c}\ \bdiamond{c}{d}{b}{d}
+\bdiamond{b}{c=b+e_\Gb}{\!\!\!\!\!\!\!\!\!d=c+e_\Ga}{c}\ \bdiamond{c}{d}{b}{d} -1
\Big] =\beta \ . }
The l.h.s. of (4) reads, if $a_{\Ga\Gb}:= (a,e_\Ga-e_\Gb)$ and 
$s(x)=\sin{\pi x\over n}$ 
\eqn\IIIk{\eqalign{ l.h.s.&= \sum_{j\ne i}
{s(a_{\Ga\Gb}+1)\over s(a_{\Ga\Gb})}
\Big[
{s(b_{\Gb\Gc}+1)\over s(b_{\Gb\Gc})}{s(c_{\Gc\Ga}+1)\over s(c_{\Gc\Ga})}
+{s(b_{\Gb\Ga}+1)\over s(b_{\Gb\Ga})}{s(c_{\Ga\Gc}+1)\over s(c_{\Ga\Gc})}-1
\Big]  \cr
&= \sum_{j\ne i}\Big[ 
{s(a_{\Ga\Gb}+1)\over s(a_{\Ga\Gb})}
{s(a_{\Gb\Gc}+1)\over s(a_{\Gb\Gc})}{s(a_{\Gc\Ga})\over s(a_{\Gc\Ga}-1)}
+{s(a_{\Ga\Gc}+2)\over s(a_{\Ga\Gc}+1)}
-{s(a_{\Ga\Gb}+1)\over s(a_{\Ga\Gb})}
\Big] \cr
&={s(a_{\Ga\Gb}+1)\over s(a_{\Ga\Gb})}
{s(a_{\Gb\Gc}+1)\over s(a_{\Gb\Gc})}{s(a_{\Gc\Ga})\over s(a_{\Gc\Ga}-1)}
+{s(a_{\Ga\Gc}+2)\over s(a_{\Ga\Gc}+1)}
-{s(a_{\Ga\Gb}+1)\over s(a_{\Ga\Gb})} \cr
&+{s(a_{\Ga\Gc}+1)\over s(a_{\Ga\Gc})}
{s(a_{\Gc\Gb}+1)\over s(a_{\Gc\Gb})}{s(a_{\Gb\Ga})\over s(a_{\Gb\Ga}-1)}
+{s(a_{\Ga\Gb}+2)\over s(a_{\Ga\Gb}+1)}
-{s(a_{\Ga\Gc}+1)\over s(a_{\Ga\Gc})}\ . \cr } }
One may check by a direct (and boring) calculation that this is
equal to $\beta$ but it is simpler to
regard it as a function of the components $a_1$ and $a_2$ of $a$; 
this is a meromorphic function; it is bounded as $a_1$, $a_2$ go
to $\infty$, and the residues at the possible poles 
($a_{\Ga\Gb}=0\ \mod n$, $a_{\Gb\Gc}=0\ \mod n$, $a_{\Ga\Gc}=-1\ \mod n$)
all vanish. By Liouville theorem, 
we conclude that this is a constant, that may be 
computed by letting $a_1$, $a_2$ go to $-i \infty$, for which we
find the result $\beta$, which establishes \IIIj\ and hence \IIIg.

We have repeated the same analysis for the case of $sl(4)$, i.e. 
$k=4$. What is quite intriguing is that the identity
\IIIk\ plays a role several times in that $k=4$ calculation: 
first to simplify 
some contributions, in the same way as \IIIha\ did, and secondly 
to help to establish the analogue of \IIIk. The relevant (very cumbersome)
trigonometric identity is again proved by making use of Liouville
theorem, and the previous identity \IIIk\ is just what guarantees the 
vanishing of the residues at the possible poles. 
There is little doubt that this structure extends to higher $k$ 
but we have not been 
able to fully disentangle the hierarchy of identities that appear there. 

Note that this approach hints at a property apparently stronger than 
the identities \IIg{}, namely that, for $k=3$ say, 
\eqn\IIIm{ (Z_4 - Z_2 G_2 - \Gb Z_1 G_3)_{ab;e} =0 }
where $(Z_L)_{ab;e}$ denotes a trace like in \IId\ but 
with two heights fixed on one side, one on the other. 
Such 
partition functions have been already the object 
of some interest in the past 
\refs{\SaBa,\Ca}.
%
%
\newsec{Conclusions and perspectives}
\nind 
In this paper we have presented some evidence that traces
of operators in representations of the Hecke algebra satisfy remarkable
identities that have a large character of universality, that is, that
are independent of the representation. Clearly a more
systematic proof is still needed. Also one should extend 
it to cases of graphs with multiple edges (as those appearing 
in orbifold graphs \Ko); we do have evidence that the identities 
remain valid in those cases, when the appropriate modifications, 
v.i.z. the introduction of additional {\it edge} degrees of freedom, are
taken into account.

Let us mention some directions in which these identities are or 
might be useful:

- As mentionned in the Introduction, 
these traces also yield contributions to the two-point function
of lattice operators, the so-called ``Pasquier operators": cf 
\refs{\VPun,\PZ}. It is
an interesting question to see if one can extend this to 3-point functions,
as it would provide some non trivial information on the operator 
algebra. This has been achieved so far only for $k=2$ \PZ, and 
there is some indication that the situation for $k\ge 3$ is more
subtle. 

- Regarding the sequences of $U$ or $g$ 
appearing in the traces \IId\ as the product along a cylinder, one 
 can also use the exact results of this paper to compute the partition
function of a quasi-one-dimensional system. We refer the interested reader
to \SL\ for a detailed discussion. 

- The original hope was that the detailed knowledge of these identities 
would provide additional and sufficient information to
fully determine the Boltzmann weights in cases  where they are 
not known. This is based on the following observation. 

A powerful method of determination of Boltzmann weights 
is through the calculation of ``cells" \refs{\Ocn,\DFZun,\Roc,\PZh}, that 
intertwine between the weights pertaining to the graph
$\CG_1$ and those relative to some basic graph $\CA$. 
Some identities, like \IIj, are such that if they are 
satisfied by the basic graphs, the mere existence 
of cells (with their orthonormality properties) guarantees
that they are also satisfied by the graph $\CG_1$. 
This is {\it not} the case with identities \IIg{}. Thus 
if these identities are satisfied both by $\CA$ and $\CG_1$, 
this implies non trivial constraints on the cells.
At this stage, the systematic exploitation of these constraints
 has not yet been achieved.
We hope to report on that matter in the future.

\vskip1cm
\noindent{\bf Acknowledgements} 
It is a pleasure to acknowledge stimulating discussions with 
P. Di Francesco, V. Pasquier and  P. Pearce.
 Y.-K. Z. is supported by the Australian
Research Council and partially by Natural Science
Foundation of China. 
J.-B. Z. wants to thank O. Foda and the Mathematics Department, 
University of Melbourne, for their generous hospitality that 
made this collaboration possible.

\appendix{A}{}
This appendix presents the proof of the assertion of eq. \IIdd, in the
case of the ``basic" graphs ${\CA}^{(n)}$. 

Let $V^{(\Gl)}$ denote the representation space of weight $\Gl$
of the algebra $\CU_q sl(k)$. We recall that for  $q$ is a root of unity, 
$q=\exp{i \pi\over n}$, it is consistent to  
 restrict to the admissible weights $\Gl\in \CP^{(n)}_{++}$. 
If $W$ is a certain representation space of the quantum group, 
we thus decompose $W$ 
as $W\cong \oplus_\Gl \ \CC_W^{(\Gl)}  \otimes V^{(\Gl)}$,
with multiplicity spaces $\CC_W^{(\Gl)} $,  and
denote $\bun^{(\Gl)}_W$ 
the projector on the irreducible representation $(\Gl)$. By Schur's lemma
any operator $x$ on $W$ that commutes with the action of the quantum group 
reads $ x=\oplus_\Gl \, x^{(\Gl)}\otimes \bun^{(\Gl)}_W$.

To describe the subspace $\CP_{\Gl_0\Gl_L}^{(L)}$ 
of paths of length $L$ on a basic
graph going from $\Gl_0$ to $\Gl_L$, it is appropriate to proceed as follows 
\item{i)} one introduces 
the space $V_0=V^{(\Gl_0)}$ and attach it to site $0$
of the chain, 
\item{ii)} one introduces $L$ copies $V_1, \cdots, V_L$ of the 
$k$-dimensional representation of $\CU_q sl(k)$, attached to 
the successive {\it edges} of the path. 

\noindent If one decomposes
\eqn\IIde{ V_0\otimes V_1 \otimes \cdots \otimes V_L 
=\oplus_\Gl \,\, \CC^{(\Gl)}  \otimes V^{(\Gl)}\ , }
then
\eqn\IIde{ \CP_{\Gl_0\Gl_L}^{(L)} =  \CC^{(\Gl_L)}  \ .}
%
Implicit in this procedure  is the connection between height
(IRF) models and ``vertex models" where the spins  take their 
values in some fixed representation of $\CU_q sl(k)$ (here the 
fundamental, $k$-dimensional one) \VPae. 

Let us denote  $W:= V_1 \otimes \cdots \otimes V_L$ for short and let 
$\Delta(\bun^{(\Gl_L)})$ be the projection on the irreducible representation 
$(\Gl_L)$ induced by the coproduct in the space $V_0\otimes W$. 
We are interested  in operators $x$ that are generated by 
the Hecke algebra and thus are in the commutant of the 
quantum group acting on $W$. Also these $x$ act as the identity on $V_0$.
For such an operator
, computing the trace in the space $V_0\otimes W$  of 
$(\bun_{0}^{(\Gl_0)}\, . x. \, \Delta(\bun^{(\Gl_L)})$
 amounts to taking the trace in the subspace 
$\CP_{\Gl_0\Gl_L}^{(L)}$ 
\eqn\IIdf{
\tr_{ V_0\otimes \cdots \otimes V_L} (\bun_{V_0}^{(\Gl_0)}\, 
. x. \, \Delta(\bun^{(\Gl_L)}) \, =
\tr_{\CP^{(L)}_{\Gl_0 \Gl_L}} (x) \, \, \tr (\bun^{(\Gl_L)}  )
\ ,}
and $\tr_{\CP^{(L)}_{\Gl_0 \Gl_L}} (x) $ is what is denoted 
$\(\tr(x)\)_{\Gl_0 \Gl_L}  $ in eq. \IIe; $\tr \bun^{(\Gl_L)}$
is the dimension of the irreducible representation $(\Gl_L)$.

By the remarks above, one may 
decompose $W=\oplus_\Gl\,\, \CC_W^{(\Gl)} \otimes V^{(\Gl)}$, 
write accordingly
 $x= {\rm id}_0 \otimes \sum_\Gl (x^{(\Gl)} \otimes \bun^{(\Gl}) $,
and compute 
\eqn\IIdi{\eqalign{
\(\tr(x)\)_{\Gl_0 \Gl_L}  \ \tr (\bun^{(\Gl_L)}) &= 
\sum_{\Gl\in \CP_{++}^{(n)}} \tr_{V_0\otimes (\CC_W^{(\Gl)}
\otimes V^{(\Gl)})} \(\bun_{0}^{(\Gl_0)} . 
( x^{(\Gl)}\otimes \bun_W^{(\Gl)}
  ) .  \Delta( \bun^{(\Gl_L)}) \) 
\cr & =
\sum_{\Gl \in \CP_{++}^{(n)}} \tr_{V_0\otimes W}  \(\bun_{V_0}^{(\Gl_0)} . 
\bun_W^{(\Gl)} . \GD(\bun^{(\Gl_L)}) \) \quad 
\tr_{\CC_W^{(\Gl)}} x^{(\Gl)} \ . \cr }}
The first trace is nothing else than the fusion coefficient 
$N_{\Gl_0 \Gl}^{\ \ \Gl_L}$ times $\tr (\bun^{(\Gl_L)})$.
If we call $x_L^{(\Gl)}$ the last trace appearing 
in \IIdi, we obtain the desired result \IIdd\ for the basic
graph, for which $N_\Gl^{(\CA)} \equiv N_\Gl$. 

\listrefs
\bye